\documentclass[runningheads]{lncse}

\usepackage{amsmath}
\usepackage{amsfonts}
\usepackage{epsfig}
\usepackage{graphics}

\begin{document}

\title{A Spectral Element Reduced Basis Method in Parametric CFD}

% your contribution title if the original one is too long use an abbreviated title (for running head):
\titlerunning{A Spectral Element Reduced Basis Method in Parametric CFD}

\author{ Martin W. Hess\inst{1} \and Gianluigi Rozza\inst{2}   }

% Use \authorrunning{Short Title} for an abbreviated version of the author list (for running head):
\authorrunning{M. Hess and G. Rozza}   

\institute{
SISSA mathLab, International School for Advanced Studies, via Bonomea 265, I-34136 Trieste, Italy, {\tt mhess@sissa.it}
\and SISSA mathLab, International School for Advanced Studies, via Bonomea 265, I-34136 Trieste, Italy, {\tt gianluigi.rozza@sissa.it}
}

\maketitle

\begin{abstract}
%Each paper should be preceded by an abstract (10--15 lines long) that summarizes the content. 
%The abstract will appear \textit{online} at {\tt www.SpringerLink.com} and be available with unrestricted access. 
%This allows unregistered users to read the abstract as a teaser for the complete paper. 
%As a general rule the abstracts will not appear in the printed version of the book unless it is the style of your particular book or that of the series to which your book belongs.

We consider the \emph{Navier-Stokes} equations in a channel with varying \emph{Reynolds} numbers.
The model is discretized with high-order spectral element \emph{ansatz} functions, resulting in $14'259$ degrees of freedom.
The steady-state snapshot solutions define a reduced order space, which allows to accurately evaluate 
the steady-state solutions for varying \emph{Reynolds} number with a reduced order model within a fixed-point iteration.
In particular, we compare different aspects of implementing the reduced order model with respect to the use of a spectral element discretization.
It is shown, how a multilevel static condensation \cite{Hess:Sherwin:2005} 
in the pressure and velocity boundary degrees of freedom can be combined with a reduced order modelling 
approach to enhance computational times in parametric many-query scenarios.

\end{abstract}

\section{Introduction}

The use of spectral element methods in computational fluid dynamics \cite{Hess:Sherwin:2005} allows highly accurate computations 
by using high-order spectral element \emph{ansatz} functions. Typically, an exponential error decay can be observed under p-refinement.
See \cite{Hess:CHQZ1}, \cite{Hess:CHQZ2}, \cite{Hess:PateraSEM}, \cite{Hess:Herrero2013132}, \cite{Hess:Taddei} for an introduction and overview of the applications.

This work is concerned with the reduced basis method (RBM, \cite{Hess:RBref}) of a channel flow, governed by the \emph{Navier-Stokes} equations, and
discretized with the spectral element method into $14'259$ degrees of freedom. In particular, we are interested in computing the steady-state solutions for varying
\emph{Reynolds} number with a reduced order model, guaranteeing competitive computational performances.

Section 2 introduces the governing equations and used fixed-point iteration algorithm. Section 3 introduces the spectral element discretization, while
section 4 describes the model reduction approach. Numerical results are provided in section 5, while section 6 summarizes and concludes the 
work by also providing new perspectives.

\section{Problem Formulation}

Let $\Omega \in \mathbb{R}^2$ be the computational domain.
Incompressible, viscous fluid motion in a spatial domain $\Omega$ over a time interval $(0, T)$
is governed by the incompressible
 \emph{Navier-Stokes} equations with velocity $\mathbf{u}$, 
pressure $p$, kinematic viscosity $\nu$ and a body forcing $f$, \eqref{Hess:NSE0} - 
\eqref{Hess:NSE1}:

\begin{eqnarray}
\frac{\partial \mathbf{u}}{\partial t} + \mathbf{u} \cdot \nabla \mathbf{u} &=& - \nabla p + \nu \Delta \mathbf{u} + f, \label{Hess:NSE0} \\
\nabla \cdot \mathbf{u} &=& 0.
\label{Hess:NSE1}
\end{eqnarray}
\noindent Boundary and initial conditions are prescribed as 
\begin{eqnarray}
\mathbf{u} &=& d \quad \text{ on } \Gamma_D \times (0, T), \\
\nabla \mathbf{u} \cdot n &=& g \quad \text{ on } \Gamma_N \times (0, T), \\
\mathbf{u} &=& u_0 \quad \text{ in } \Omega \times 0,
\label{Hess:NSE_boundaryCond}
\end{eqnarray}

\noindent with $d$, $g$ and $u_0$ given and $\partial \Omega = \Gamma_D \cup \Gamma_N$, $\Gamma_D \cap \Gamma_N = \emptyset$.
The \emph{Reynolds} number $Re$ depends on the viscosity $\nu$ through the characteristic velocity $U$ and characteristic length $L$ via 
$Re = \frac{UL}{\nu}$, \cite{Hess:Pitton2017534}.

%\begin{eqnarray}
% Re = \frac{UL}{\nu}.
% \label{Hess:Re_visc}
%\end{eqnarray}

In particular, we are interested in computing 
the steady states for varying viscosity $\nu$, 
such that $\frac{\partial \mathbf{u}}{\partial t} = 0$.
A solution $\mathbf{u}(\nu_1)$ for a parameter value $\nu_1$, 
can be used as an initial guess for a fixed point iteration 
to obtain the steady state solution $\mathbf{u}(\nu_2)$ at a 
parameter value $\nu_2$, provided that the solution 
$\mathbf{u}(\nu)$ depends continuously 
on $\nu$ in the interval $[\nu_1, \nu_2]$.

%momentum and pressure equation....

\subsection{Oseen-Iteration}

The \emph{Oseen}-iteration is a secant modulus fixed-point iteration, which in general exhibits a linear rate of convergence \cite{Hess:Oseen}.
Given a current iterate (or initial condition) $\mathbf{u}^k$, the linear system

%Given $u^k \in V$, find $u^{k+1}$, such that
%\begin{eqnarray}
%  a(u^{k+1}, v) + b(u^k, u^{k+1}, v) = (f,v) \quad \forall v \in V
%\end{eqnarray}

\begin{eqnarray}
 -\nu \Delta \mathbf{u} + (\mathbf{u}^k \cdot \nabla) \mathbf{u} + \nabla p &=& f  \text{ in } \Omega, \label{Hess:eq_Oseen_main}\\
\nabla \cdot \mathbf{u} &=& 0   \text{ in } \Omega, \\
 %\mathbf{u} &=& f_d  \text{ on } \partial \Omega,
\mathbf{u} &=& d \quad \text{ on } \Gamma_D, \\
\nabla \mathbf{u} \cdot n &=& g \quad \text{ on } \Gamma_N, 
\end{eqnarray}
\noindent is solved for the next iterate $\mathbf{u}^{k+1} = \mathbf{u}$. 
A typical stopping criterion is that the relative change between iterates in the $H^1$ norm falls below a predefined tolerance.
An initial solution $\mathbf{u}^0(\nu_0)$ is computed by time-advancement of \eqref{Hess:NSE0}--\eqref{Hess:NSE1} from zero initial conditions at a parameter value $\nu_0$,
and the whole parameter domain is then explored by using a continuation method with the \emph{Oseen}-iteration.

%The \emph{Oseen}-iteration is then initialized with $\mathbf{u}(\nu_0)$ to solve for the velocity at a value close to $\nu_0$. 
%Solutions over the whole parameter domain are then explored by using a continuation method, i.e., using 
%the solution of the fixed point iteration for the previous parameter value as initial condition at the next parameter value.
%This improves computation time over initializing the \emph{Oseen}-iteration with a zero field.

\section{Spectral Element Discretization}

The \emph{Navier-Stokes} problem is discretized with the spectral element 
method. The spectral/hp element software framework used is Nektar++ in version 4.3.5, \cite{Hess:nektar}\footnote{See \textbf{www.nektar.info}.}.
%An initial solution is computed with a time-dependent simulation.
%This initial solution then serves as an initial guess for the Oseen-iteration. 
The discretized system to solve in each step of the \emph{Oseen}-iteration is given by \eqref{Hess:fully_expanded} as
\begin{eqnarray}
\begin{bmatrix}
\begin{array}{ccc}
A & -D^T_{bnd} & B  \\
-D_{bnd} & 0 & -D_{int} \\
\tilde{B}^T & -D^T_{int} & C
\end{array}
\end{bmatrix}
\begin{bmatrix}
\begin{array}{ccc}
v_{bnd} \\
p \\
v_{int} 
\end{array}
\end{bmatrix}
&=
\begin{bmatrix}
\begin{array}{ccc}
f_{bnd} \\
0 \\
f_{int}
\end{array}
\end{bmatrix},
\label{Hess:fully_expanded}
\end{eqnarray}
\noindent where $v_{bnd}$ and $v_{int}$ denote velocity degrees of freedom on the boundary and in the interior, respectively.
Correspondingly, $f_{bnd}$ and $f_{int}$ denote forcing terms on the boundary and interior, respectively.
The matrix $A$ assembles the boundary-boundary coupling, $B$ the boundary-interior coupling, $\tilde{B}$ the interior-boundary coupling and
$C$ assembles the interior-interior coupling of elemental velocity \emph{ansatz} functions. 
In the case of a \emph{Stokes} system, it holds that $B = \tilde{B}^T$, but this is not the case for the \emph{Oseen} equation, since 
the linearization term $(\mathbf{u}^k \cdot \nabla) \mathbf{u}$ is present in \eqref{Hess:eq_Oseen_main}.
The matrices $D_{bnd}$ and $D_{int}$ assemble the pressure-velocity boundary and pressure-velocity interior contributions, respectively.

The linear system \eqref{Hess:fully_expanded} is assembled in local degrees of freedom, resulting in block matrices $A, B, \tilde{B}, C, D_{bnd}$ and $D_{int}$,
each block corresponding to a spectral element. In particular, this means that the system is singular in this form.
To solve the system, the local degrees of freedom need to be gathered into the global degrees of freedom \cite{Hess:Sherwin:2005}.
Since $C$ contains the interior-interior contributions, it is invertible and the system can be statically condensed into
\begin{eqnarray}
\begin{bmatrix}
\begin{array}{ccc}
A - B C^{-1} \tilde{B}^T & B C^{-1} D^T_{int}-D^T_{bnd} & 0  \\
D_{int} C^{-1} \tilde{B}^T-D_{bnd} & -D_{int} C^{-1} D^T_{int} & 0 \\
\tilde{B}^T & -D^T_{int} & C
\end{array}
\end{bmatrix}
&\begin{bmatrix}
\begin{array}{ccc}
v_{bnd} \\
p \\
v_{int} 
\end{array}
\end{bmatrix} \nonumber \\
=&
\begin{bmatrix}
\begin{array}{ccc}
f_{bnd} - B C^{-1} f_{int} \\
D_{int} C^{-1} f_{int} \\
f_{int}
\end{array}
\end{bmatrix}.
\label{Hess:first_condensed}
\end{eqnarray}
By taking the top left $2 \times 2$ block and reordering the degrees of freedom such that the mean pressure mode of each element is inserted into the 
corresponding block of $\hat{A}$ results in
\begin{eqnarray}
\begin{bmatrix}
\begin{array}{cc}
\hat{A} & \hat{B}  \\
\hat{C} & \hat{D}
\end{array}
\end{bmatrix}
\begin{bmatrix}
\begin{array}{ccc}
b \\
\hat{p}
\end{array}
\end{bmatrix}
&=
\begin{bmatrix}
\begin{array}{ccc}
\hat{f}_{bnd} \\
\hat{f}_p
\end{array}
\end{bmatrix},
\label{Hess:hats}
\end{eqnarray}
\noindent where $\hat{D}$ is invertible, such that a second level of static condensation can be employed.
We have:
\begin{eqnarray}
\begin{bmatrix}
\begin{array}{cc}
\hat{A} - \hat{B} \hat{D}^{-1} \hat{C} & 0  \\
\hat{C} & \hat{D}
\end{array}
\end{bmatrix}
\begin{bmatrix}
\begin{array}{ccc}
b \\
\hat{p}
\end{array}
\end{bmatrix}
&=
\begin{bmatrix}
\begin{array}{ccc}
\hat{f}_{bnd} - \hat{B} \hat{D}^{-1} \hat{f}_p \\
\hat{f}_p
\end{array}
\end{bmatrix}.
\label{Hess:hats_final}
\end{eqnarray}
When the vector $b$ is computed, which contains the velocity boundary degrees of freedom and the mean pressure modes,
the remaining solution components are computed by reverting the steps of the static condensations \cite{Hess:Sherwin:2005}.
The main computational effort lies in solving the final system \eqref{Hess:hats_very_final}
\begin{eqnarray}
(\hat{A} - \hat{B} \hat{D}^{-1} \hat{C})  b  = \hat{f}_{bnd} - \hat{B} \hat{D}^{-1} \hat{f}_p .
\label{Hess:hats_very_final}
\end{eqnarray}
Additionally, the matrices $C$ and $\hat{D}$ need to be inverted, which due to the elemental block structure requires 
inverting submatrices in the size of the degrees of freedom per element for each submatrix.

%Then look at the form (17) and discuss

%write section 10.1.3 from Nektar user guide

\section{Reduced Order Modelling}

The reduced order model (ROM) aims to represent the full order solution accurately in the parameter domain of interest.
Two ingredients are essential to RB modelling, a projection onto a low order space of snapshot solutions and 
an offline-online decomposition for computational efficiency, \cite{Hess:Lassila2014}.
A set of snapshots is generated by solving \eqref{Hess:hats_very_final} over a coarse discretization of the parameter domain 
and used to define a projection space $U$ of size $N$. 
The proper orthogonal decomposition computes a singular value decomposition of the snapshot solutions to $99.9\%$ of the most dominant modes \cite{Hess:RBref}, which
defines the projection matrix $U \in \mathbb{R}^{N_\delta \times N}$ to project system \eqref{Hess:hats_very_final} and obtain 
the reduced order solution $b_N$.

%such that the reduced order system can be expressed as

%\begin{eqnarray}
%U^T (\hat{A} - \hat{B} \hat{D}^{-1} \hat{C}) U  b_N  = U^T (\hat{f}_{bnd} - \hat{B} \hat{D}^{-1} \hat{f}_p) .
%\label{Hess:hats_very_final_red}
%\end{eqnarray}

%\noindent with the reduced order solution $b_N$.

%While the twice statically condensed system has the lowest dimension, it is difficult to explicitly extract the parameter-dependence from this
%system, which is present in every matrix and required for the offline-online decomposition.
%Instead the offline-online decomposition will proceed from one level of static condensation, i.e., \eqref{Hess:first_condensed}.

\subsection{Offline-Online Decomposition}

The offline-online decomposition \cite{Hess:RBref} allows fast input-output evaluations independent of the original model size $N_\delta$.
It is a crucial part of an efficient reduced order model but since the static condensation includes 
the inversion of the parameter-dependent matrix $C$, an intermediate projection is introduced.

The reduced order model considers the top left $2 \times 2$ block of \eqref{Hess:first_condensed}, i.e., one level of static
condensation \cite{Hess:Sherwin:2005}. 
During the offline phase, full-order solutions have been computed over the parameter domain of interest, which now serve as a projection space 
to define the reduced order setting. This projection space incorporates the transformation of local velocity boundary degrees of freedom 
to global velocity boundary degrees of freedom and the reordering of mean pressure degrees of freedom. 
The projection space then takes the form $U = PMV$ with a permutation matrix $P$ to reorder the degrees of freedom and a transformation $M$ from local to global degrees of freedom.
The collected offline data $V$ contain the gathered velocity and mean pressure modes as well as interior pressure modes.

%\begin{eqnarray}
% U = P
%\begin{bmatrix}
%\begin{array}{cc}
%M & 0 \\
%0 & Id
%\end{array}
%\end{bmatrix}
%V,
%\end{eqnarray}

%\noindent with a permutation matrix $P$ to reorder the degrees of freedom and a transformation $M$ from local to global degrees of freedom.
%The collected offline data $V$ contain the gathered velocity and mean pressure modes as well as interior pressure modes.

The projected system is then of the form 
\begin{eqnarray}
A_N = U^T
 \begin{bmatrix}
\begin{array}{cc}
A - B C^{-1} \tilde{B}^T & B C^{-1} D^T_{int}-D^T_{bnd}  \\
D_{int} C^{-1} \tilde{B}^T-D_{bnd} & -D_{int} C^{-1} D^T_{int} 
\end{array}
\end{bmatrix}
U,
\label{Hess:red_general}
\end{eqnarray}
\noindent and upon its solution, the interior velocity dofs can be computed by resubstituting into \eqref{Hess:first_condensed} at the reduced order level.
To achieve fast reduced order solves, the offline-online decomposition expands \eqref{Hess:red_general} in the parameter of interest and 
computes the parameter independent projections offline to be stored as small-sized matrices of the order $N \times N$.
Since in an \emph{Oseen}-iteration each matrix is dependent on the previous iterate, the submatrices corresponding to each basis function is assembled and then 
formed online using the reduced basis coordinate representation of the current iterate. 
This is analogous to reduced order assembly of the nonlinear term in the \emph{Navier-Stokes} case, \cite{Hess:Lassila2014}.

%This is made explicit here on the example of the matrix $A$ and analogous for the other cases.
%Since $A$ denotes the boundary-boundary coupling of velocity degrees of freedom, each element $(i,j)$ of $A$ is of the form

%\begin{eqnarray}
% A_{ij} = (\nabla \Phi^b_i, \nu \nabla \Phi^b_i) + (\nabla \Phi^b_i, u_k \cdot \nabla \Phi^b_j),
%\end{eqnarray}

%\noindent with spectral element \emph{ansatz} function on the boundary $\Phi^b_i$ and current \emph{Oseen} iterate $u_k$.
%Since the viscosity is considered a parameter, $\nu$ is separated and the parameter-independent part $(\nabla \Phi^b_i, \nabla \Phi^b_i)$ is 
%assembled offline.

%introduce reduced quantities here

%snapshot space generation

%\subsection{Reduced order Model within static condensation}

%A further gain can be achieved by not only using the ROM to speed up the computation time for the linear system \eqref{Hess:hats_final}, but
%also when inverting the block matrices in $D$ and $\hat{D}$.

%A direct approach did not yield good results...

%instead lets do this special thingy here....

\section{Model and Numerical Results}

\begin{figure}[t]
\includegraphics[scale=.35]{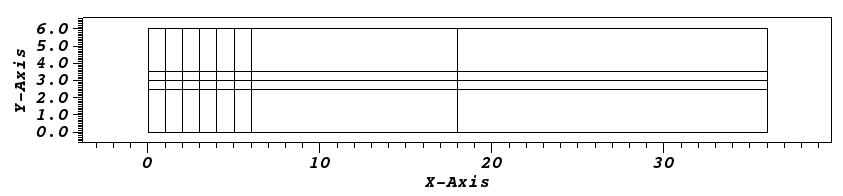}
\caption{Rectangular domain for the channel flow, shown is the $4 \times 8$ spectral element grid. }
\label{Hess:domain}       % Give a unique label
\end{figure}

We consider a channel flow in the domain considered in Fig.~\ref{Hess:domain}, similar to the model considered in \cite{Hess:Pitton2017}.
The rectangular domain $\Omega(x,y) = [0,36] \times [0,6]$ is decomposed into $32$ spectral elements. 
The spectral element expansion uses modal Legendre polynomials of order $p = 12$ in the velocity. 
The pressure \emph{ansatz} space is chosen of order $p-2$ to fulfill the inf-sup stability condition, \cite{Hess:BBF}, \cite{Hess:infsup}.
The inflow is defined for $y \in [2.5, 3.5]$ as $u_x(0,y) = (y-2.5)(3.5-y)$.
At $x = 36$ is the outflow boundary, everywhere else are zero velocity walls.
Note that the velocity boundary degrees of freedom are along the boundaries of the spectral elements and not only the domain boundary, resulting
in $3072$ local degrees of freedom for this problem.
This is a simplified model of a contraction-expansion channel \cite{Hess:Pitton2017534}, where flow occurs though a narrowing of variable width.
Variations in the width have been moved to variations in the \emph{Reynolds} number and only the section after the narrowing comprises the 
computational domain.
The relation to the \emph{Reynolds} number %\eqref{Hess:Re_visc} 
is established with $U=\frac{1}{4}$ as the maximum inflow velocity and 
$L=1$ as the width of the narrowing as $Re = \frac{1}{4\nu}$.

%\begin{eqnarray}
% Re = \frac{1}{4\nu}.
% \label{Hess:Re_visc_example}
%\end{eqnarray}

\begin{figure}[ht]
 \includegraphics[scale=.15]{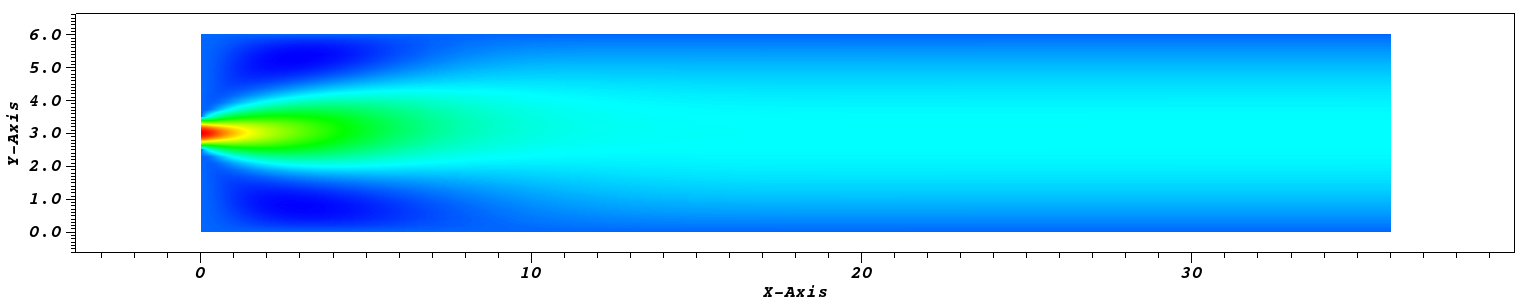} $\quad$
 \includegraphics[scale=.2]{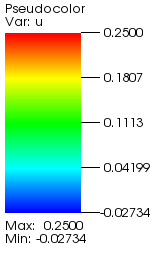} \\
 \includegraphics[scale=.15]{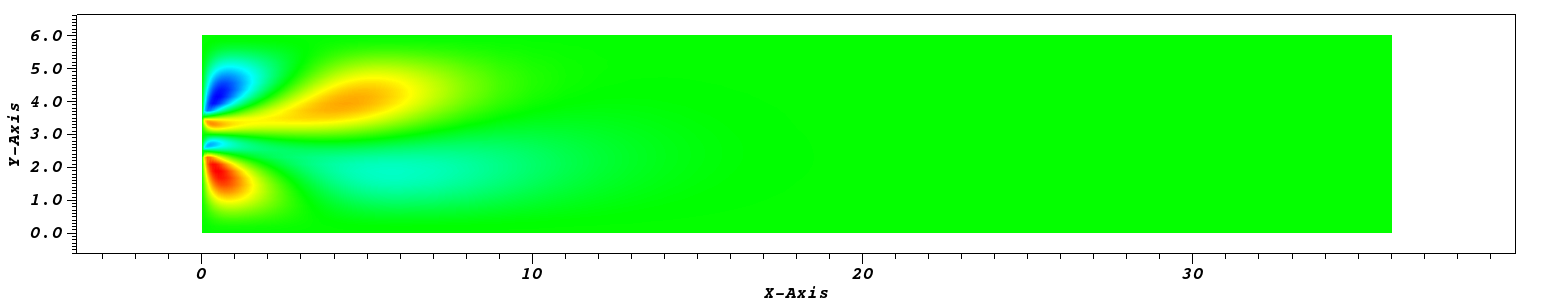} $\quad$
 \includegraphics[scale=.2]{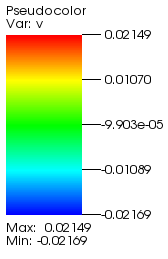}
 \caption{Full order, steady-state solution for $\nu = 0.0075$. Top is the velocity in x-direction, below is the velocity in y-direction.}
 \label{Hess:Xp0075}
\end{figure}

\begin{figure}[ht]
 \includegraphics[scale=.15]{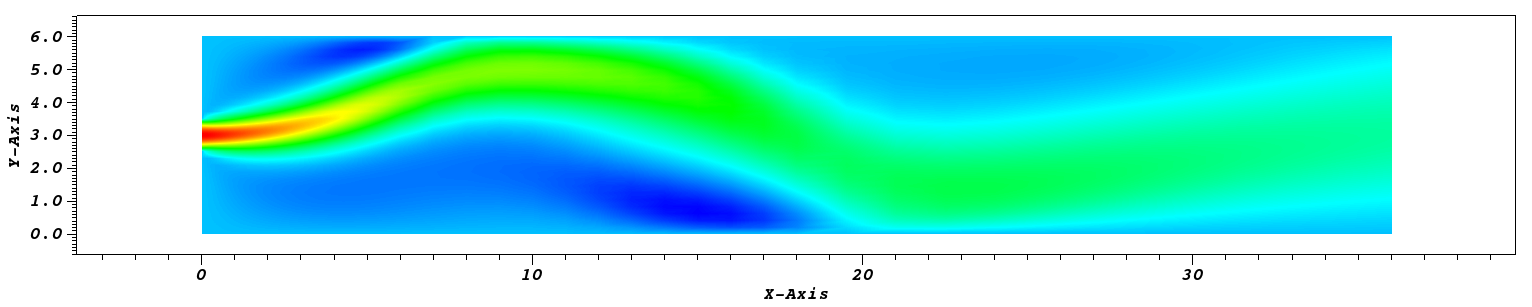} $\quad$
 \includegraphics[scale=.2]{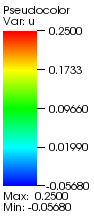} \\
 \includegraphics[scale=.15]{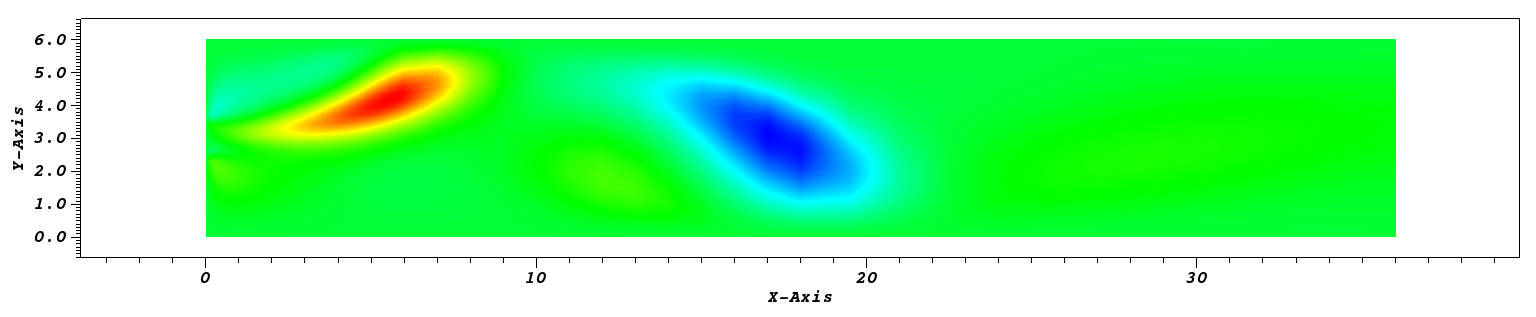} $\quad$
 \includegraphics[scale=.2]{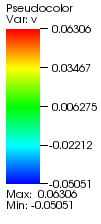}
 \caption{Full order, steady-state solution for $\nu = 0.0025$. Top is the velocity in x-direction, below is the velocity in y-direction.}
 \label{Hess:Xp0025}
\end{figure}

\begin{figure}[ht]
 \center
 \includegraphics[scale=.9]{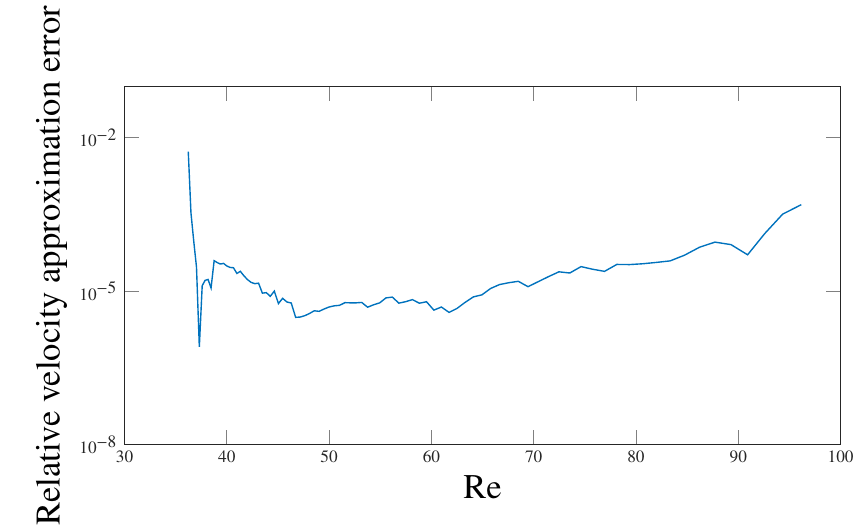}
 \caption{Relative error in the velocity over the parameter domain.}
 \label{Hess:relerr}
\end{figure}

Consider a parametric variation in the viscosity $\nu$, ranging from $\nu = 0.0075$ to $\nu = 0.0025$, which corresponds 
to \emph{Reynolds} numbers between $33$ and $100$.
The solution for  $\nu = 0.0075$ is shown in Fig.~\ref{Hess:Xp0075}.
It is slightly unsymmetrical, which marks the onset of the \emph{Coanda} effect \cite{Hess:wille_fernholz_1965}, \cite{Hess:Coanda}, which is a known phenomenon characterized as a `wall-hugging' effect
occurring at these \emph{Reynolds} numbers. 
The solution for  $\nu = 0.0025$ is shown in Fig.~\ref{Hess:Xp0025}.
Here, the \emph{Coanda} effect is fully developed as the flow orients itself along the boundaries.
Using model reduction with the form \eqref{Hess:first_condensed}, which allows the offline-online decomposition 
or using form \eqref{Hess:hats_very_final}, which has the lowest full-order system size, resulted in similar computational results.
Shown in Fig.~\ref{Hess:relerr} is the relative $H^1_0(\Omega)$ error in the velocity between the full order and reduced order model.

While the full-order solves were computed with Nektar++, the reduced-order computations were done in a separate python code.
To compare computational gains, compute times between a full order solve and a reduced order solve both implemented in python 
are taken. The compute times reduce by a factor of $50$, i.e., for a single iteration step from about $40$s to under $1$s.
Current work also aims to extend the software to make it available as a SEM-ROM software framework
within the AROMA-CFD project (see Acknowledgment) as ITHACA-SEM.

%model description with picture and one snapshot -- maybe later

%do the comparison

%expect two graphs 
%- fast approx with standard large RB
%- faster approx with 'statically condensed RB'
%  maybe additional cost of element wise C-solves

\section{Conclusion and Outlook}

It has been shown that the reduced basis technique generates accurate reduced order models of small size for 
channel flow discretized with spectral elements up to a \emph{Reynolds} number of $100$.
The use of basis functions obtained by the spectral element method suggests a potential important synergy between
high-order and reduced basis methods, see also \cite{Hess:Taddei}.
Due to the multilevel static condensation used here, particular care must be taken to achieve an offline-online decomposition.
The domain decomposition into spectral elements shows a resemblance to reduced basis element methods (RBEM),\cite{Hess:Maday2002195} \cite{Hess:lovgren_maday_ronquist_2006}.
A comparison of both approaches could be the subject of further investigation.

%when possible to statically condensed RB
%LINK with RBEM ??

%\section*{References} %
%Please do not use a bibitex file but use the reference style of the present bibliography.
%For BibTeX users, please use the style {\tt vmams} (see bibliography section) and copy and paste the generated bbl file.

%Remember that all references should be cited in the text 
%\cite{author_mini4:clar:eke:2,author_mini4:clar:eke,author_mini4:Francis:Vismath95,author_mini4:helman:RDV,author_mini4:helman:VVF} and also 
%\cite{author_mini4:mich:tar,author_mini4:rab,author_mini4:tar}.

%
%\section*{Appendix}
%\addcontentsline{toc}{section}{Appendix}
%
%
%When placed at the end of a chapter or contribution (as opposed to at the end of the book), the numbering of tables, figures, and equations in the appendix section continues on from that in the main text. Hence please \textit{do not} use the \verb|appendix| command when writing an appendix at the end of your chapter or contribution. If there is only one the appendix is designated ``Appendix'', or ``Appendix 1'', or ``Appendix 2'', etc. if there is more than one.

%\begin{equation}
%a \times b = c.
%\end{equation}

\def\cprime{$'$}

\section*{Acknowledgments}  

This work was supported by European Union Funding for Research and Innovation through the European Research Council
(project H2020 ERC CoG 2015 AROMA-CFD project 681447, P.I. Prof. G. Rozza).

% 1. BibTeX version: 
%
% If you want to choose this version uncomment the next two lines and insert the name of the BibTeX file.
%
%
% \bibliographystyle{vmams}  
% \bibliography{example}        % use BibTeX file example.bib (BibTeX -> author.bbl)

\begin{thebibliography}{99}
\parskip1.0ex

\bibitem{Hess:Sherwin:2005}
{\sc G. Karniadakis and S. Sherwin},
{\em Spectral/hp Element Methods for Computational Fluid Dynamics},
Oxford University Press, 2nd ed. (2005).


\bibitem{Hess:CHQZ1}
{\sc C. Canuto, M.Y. Hussaini, A. Quarteroni and Th.A. Zhang}
{\em Spectral {M}ethods {F}undamentals in {S}ingle {D}omains},
Springer -- Scientific {C}omputation, (2006).


\bibitem{Hess:CHQZ2}
{\sc C. Canuto, M.Y. Hussaini, A. Quarteroni and Th.A. Zhang},
{\em Spectral {M}ethods {E}volution to {C}omplex {G}eometries and {A}pplications to {F}luid {D}ynamics},
Springer -- Scientific {C}omputation, (2007).


\bibitem{Hess:PateraSEM}
{\sc A.T. Patera},
{\em A Spectral Element Method for Fluid Dynamics; Laminar Flow in a Channel Expansion},
Journal of Computational Physics, {\bf 54}:3 (1984), 468--488.

\bibitem{Hess:Herrero2013132}
{\sc H. Herrero and Y. Maday and F. Pla},
{\em RB (Reduced Basis) for RB (Rayleigh--B\'enard)},
Computer Methods in Applied Mechanics and Engineering, {\bf 261--262},
(2013), 132--141.

\bibitem{Hess:Taddei}
{\sc L. Fick, Y. Maday, A. Patera and T. Taddei},
{\em A  Reduced  Basis Technique for Long-Time Unsteady Turbulent Flows},
Journal of Computational Physics (submitted), 
arXiv:https://arxiv.org/pdf/1710.03569.pdf.


\bibitem{Hess:RBref}
{\sc J.S. Hesthaven, G. Rozza and B. Stamm},
{\em Certified Reduced Basis Methods for Parametrized Partial Differential Equations},
SpringerBriefs in Mathematics, (2016).

\bibitem{Hess:Oseen}
{\sc Martin Burger},
{\em Numerical Methods for Incompressible Flow, Lecture Notes},
UCLA (2010).




%\bibitem{Hess:dijkstra_wubs_cliffe_doedel_dragomirescu_eckhardt_gelfgat_hazel_lucarini_salinger_etal._2014}
%{\sc Dijkstra, Henk A. and Wubs, Fred W. and Cliffe, Andrew K. and Doedel, Eusebius and Dragomirescu, Ioana F. and Eckhardt, Bruno and Gelfgat, Alexander Yu. and Hazel, Andrew L. and Lucarini, Valerio and Salinger, Andy G. and et al.}, 
%{\em Numerical Bifurcation Methods and their Application to Fluid Dynamics: Analysis beyond Simulation}, 
%Communications in Computational Physics, {\bf 15}: 1, 
%(2014), 1--45.

\bibitem{Hess:nektar}
{\sc C.D. Cantwell, D. Moxey, A. Comerford, A. Bolis, G. Rocco, G. Mengaldo, D. de Grazia, S. Yakovlev, J.-E. Lombard, D. Ekelschot, B. Jordi, H. Xu, Y. Mohamied, C. Eskilsson, B. Nelson, P. Vos, C. Biotto, R.M. Kirby, S.J. Sherwin},
{\em Nektar++: An open-source spectral/hp element framework},
Computer Physics Communications, {\bf 192},
(2015), 205--219.





%\bibitem{Hess:podcvtlang}
%{\sc S. Ullmann and J. Lang},
%{\em POD and CVT Galerkin reduced-order modelling of the flow around a cylinder},
%PAMM, (2012), {\bf 12}.

\bibitem{Hess:BBF}
{\sc D. Boffi, F. Brezzi and M. Fortin},
{\em Mixed Finite Element Methods and Applications},
Springer Series in Computational Mathematics, (2013).


\bibitem{Hess:infsup}
{\sc A. Quarteroni, A. Valli},
{\em Numerical  Approximation  of  Partial  Differential  Equations},
Springer-Verlag, Berlin-Heidelberg, (1994).


\bibitem{Hess:Pitton2017534}
{\sc G. Pitton, A. Quaini and G. Rozza},
{\em Computational Reduction Strategies for the Detection of Steady Bifurcations in Incompressible Fluid-Dynamics: Applications to \emph{Coanda} Effect in Cardiology},
Journal of Computational Physics, {\bf 344}, (2017), 534--557.


\bibitem{Hess:Pitton2017}
{\sc G. Pitton and G. Rozza},
{\em On the Application of Reduced Basis Methods to Bifurcation Problems in Incompressible Fluid Dynamics},
Journal of Scientific Computing, (2017).

\bibitem{Hess:wille_fernholz_1965}
{\sc R.~Wille and H.~Fernholz}, 
{\em Report on the first European Mechanics Colloquium, on the \emph{Coanda} effect}, 
Journal of Fluid Mechanics, {\bf 23}:4 
(1965), 801--819.

\bibitem{Hess:Coanda}
{\sc A. Quaini, R. Glowinski, S. \v Cani\'c},
{\em A  computational  study  on  the  generation  of  the  \emph{Coanda} effect in a mock heart chamber},
RIMS K\^oky\^uroku series, No. 2009-4, (2016).






%\bibitem{Hess:10.1080/10618562.2016.1144877}
%{\sc A. Quaini, R. Glowinski and S. Čanić},
%{\em Symmetry Breaking and Preliminary Results about a Hopf Bifurcation for Incompressible Viscous Flow in an Expansion Channel},
%International Journal of Computational Fluid Dynamics, {\bf 30}:1
%(2016), 7--19.





\bibitem{Hess:Lassila2014}
{\sc T.~Lassila, A.~Manzoni, A.~Quarteroni and G.~Rozza}, 
{\em Model Order Reduction in Fluid Dynamics: Challenges and Perspectives},
Reduced Order Methods for Modelling and Computational Reduction, Springer International Publishing, MS\&A, Vol. 9, A. Quarteroni, G.Rozza eds. 
 (2014), 235--273.



%\bibitem{Hess:Yano20130036}
%{\sc M.~Yano and A.~T.~Patera}, 
%{\em A Space{\textendash}Time Variational Approach to Hydrodynamic Stability Theory},
%Proceedings of the Royal Society of London A: Mathematical, Physical and Engineering Sciences {\bf 469}:2155
% (2013).


\bibitem{Hess:Maday2002195}
{\sc Y.~Maday and E.~M.~Ronquist}, 
{\em A Reduced-Basis Element Method},
Comptes Rendus Mathematique {\bf 335}:2
 (2002), 195--200.


\bibitem{Hess:lovgren_maday_ronquist_2006}
{\sc A.E.~Lovgren, Y.~Maday and E. M. Ronquist}, 
{\em A Reduced Basis Element Method for the Steady Stokes Problem},
ESAIM: Mathematical Modelling and Numerical Analysis {\bf 40}:3
 (2006), 529--552.


\end{thebibliography}

\ifx\undefined\bysame
\newcommand{\bysame}{\leavevmode\hbox to3em{\hrulefill}\,}
\fi

\end{document}